\documentclass[12 pt]{report}
\usepackage{amssymb,latexsym,epsfig,mathrsfs,graphpap,graphics,amsmath,amssymb}
\usepackage{amsmath}

\begin{document}

\vspace{.2in}\parindent=0mm

\begin{flushleft}
  {\bf\Large { Novel Special Affine Wavelet Transform  \vspace{.1in}  and   Associated  Uncertainity Inequalities }}
\end{flushleft}

\parindent=0mm \vspace{.2in}

  {\bf{Owais Ahmad$^{1  }$ and  Neyaz A. Sheikh$^{2}$}}

\parindent=0mm \vspace{.2in}
{{\it $^{1}$Department of  Mathematics,  National Institute of Technology, Hazratbal, Srinagar -190 006, Jammu and Kashmir, India, E-mail: $\text{siawoahmad@gmail.com}$}}

\parindent=0mm \vspace{.1in}
{{\it $^{2}$Department of  Mathematics,  National Institute of Technology, Hazratbal, Srinagar -190 006, Jammu and Kashmir, India, E-mail: $\text{neyaznit@yahoo.co.in}$}}

\parindent=0mm \vspace{.2in}
{\small {\bf Abstract.} Due to the extra degrees of freedom, special affine Fourier transform (SAFT) has achieved a respectable status within a short span and got versatile applicability in the areas of signal processing, image processing,sampling theory, quantum mechanics. However, due to its global kernel, SAFT  fails to obtain local information of non-transient signals. To overcome this, we  in this paper  introduce the concept of  novel special affine wavelet transform  (NSAWT) and extend key harmonic  analysis results to NSAWT analogous to those for the wavelet transform. We first establish some fundamental properties including Moyal's principle, Inversion formula and the range theorem.
Some Heisenberg type inequalities and Pitt's inequality are established for SAFT and consequently Heisenberg uncertainity principle is derived for NSAWT.

\parindent=0mm \vspace{.1in}
{\bf{Keywords:}} Special affine Fourier transform; Chirp Modulation; Wavelet transform;  Uncertainty principle; Convolution,

\parindent=0mm \vspace{.1in}
{\bf {Mathematics Subject Classification:}}  42C40; 42B10; 65R10; 42C15}

\parindent=0mm \vspace{.1in}
{\bf {1. Introduction}}

\parindent=0mm \vspace{.1in}
The special affine wavelet transform (SAFT), which was introduced in \cite{AS}, is a six-parameter class of linear integral transformation which generalizws several well known unitary transformations including  the Fourier transform, the fractional Fourier transform, the linear canonical transform (LCT)  and  the Fresnel transform \cite{Al,c7,JA}. The SAFT can be regarded as a time-shifted and frequency modulated version of the well known linear canonical transform \cite{cc,b1}.  Let $f^*$ denote the complex-conjugate of $f$ and $\langle f,g\rangle = \int f(x)g^*(x) dx$ be the standard $L^2$ inner product. The SAFT is a mapping ${\mathcal F}_{SAFT}: f \rightarrow {\widehat{f}}_{\bf\Lambda_S}$ and is defined as \\ 
\begin{align*}
{\mathcal F}_{SAFT}\big[f\big](\omega)=\widehat f_{\bf\Lambda_S}(\omega)=\begin{cases}\left\langle f(x)\,\mathcal{K}_{\bf \Lambda_S}(x,\omega)\,\right\rangle, & B\neq 0 \\\\
 \sqrt{D}\exp\left\{\frac{i}{2}\big(CD(\omega-p)^2+2\omega q\big)\right\}f\big(D(\omega-p)\big), & B=0,\end{cases}\tag{1.1}
\end{align*}

\parindent=0mm \vspace{.0in}
 where ${\mathcal K}_{\bf\Lambda_S}\big(x,\omega\big)$ denotes the kernel of the SAFT given by
{\small
\begin{align*}
{\mathcal K}_{\bf\Lambda_S}\big(x,\omega\big)=K_{B}^*\exp\left\{\dfrac{i}{2B}\Big(Ax^2+2x(p-\omega)-2\omega(Dp-Bq)+D\omega^2\big)\Big)\right\},~ K_{B}=\dfrac{1}{\sqrt{2\pi B}}.\tag{1.2}
\end{align*}}
and ${\bf\Lambda_S}^{(2 \times 3)}$  denotes the augmented SAFT parameter matrix, which is of the form
$${\bf\Lambda_S}^{(2 \times 3)}=\big[\bf \Lambda|~\lambda\big],$$
 which in turn is obtained by LCT matrix  $\bf\Lambda=\left[\begin{array}{cc}A & B \\C & D \\\end{array}\right]$  and an  offset vector $\lambda=\left[\begin{array}{cc}p \\q\\\end{array}\right]$ .
 
\parindent=0mm \vspace{.0in}
This is the reason that the SAFT is also called as the offset linear canonical transform. Moreover, we shall only consider the case $B\ne 0$, since the SAFT is just a chirp multiplication operation in case $B=0$. We also note that the phase-space transform (1.1) is lossless if and only if the matrix $\bf \Lambda$ is unimodular, that is; $AD-BC=1$ and for this reason, SAFT is also known as the inhomogeneous canonical transform \cite{AS1}. By virtue of the additive propery of SAFT, the inverse SAFT corresponding to (1.1) is defined by
\begin{align*}
f(x)=\langle\widehat f_{\bf\Lambda_S}(\omega)\,\mathcal{K}_{\bf \Lambda_S^{inv}}(\omega,x)\rangle\tag{1.3}
\end{align*}

\parindent=0mm \vspace{.0in}

where 
$${\bf\Lambda_S^{inv}}=\left[\begin{array}{ccc}\phantom{-}D &~ -B :& Bq-Dp\\-C & \phantom{-}A  : &Cp-Aq \\ \end{array}\right]$$The Parseval's formula for the special affine Fourier transform reads as follows
\begin{align*}
\big\langle f(x),\, g(x)\big\rangle_{L^2(\mathbb R)}=\left\langle \widehat f_{\Lambda_S}(\omega),\, \widehat g_{\Lambda_S}(\omega)\right\rangle_{L^2(\mathbb R)},\qquad \forall~ f,g\in L^2(\mathbb R).
\end{align*}

\parindent=8mm \vspace{.0in}
The theory of  wavelet transforms have emanated as a broadly used tool  in various disciplines of science and engineering including    image processing, spectrometry,machine learning,  turbulence, computer graphics, telecommunications, DNA sequence analysis,  quantum physics,  solution of  differential equations. For any $f\in L^2(\mathbb R)$, the continuous wavelet transform (CWT) is denoted by ${\mathcal W}_{\psi}[f](t,\zeta)$ and is defined as
\begin{align*}
{\mathcal W}_{\psi}\big[f\big](t,\zeta)=\dfrac{1}{\sqrt \zeta}\int_{\mathbb R}f(x)\,\overline{\psi\left(\dfrac{x-t}{\zeta}\right)}\,dx,\quad t\in\mathbb R,\, \zeta\in\mathbb R^+,\tag{1.4}
\end{align*}
where $\zeta$ is the scaling parameter and $t$ is the translation parameter.  Shah {\it et.al} \cite{shah}  introduced an amalgam of CWT and SAFT namely special affine wavelet transform (SAWT) which provides a joint time and frequency localization of signals. Covolution plays a pivotal role as far as applications of integral transforms are concerned. SAFT does not work well with the standard convolution operation. Xiang and Qin \cite{b2} introduced a convolution which works well for the SAFT and by which the SAFT of the convolutionof two functions is the product of their SAFT's and a phase factor but their convolution structure does not work well with the inverse transform. Bhandari and Zayed \cite{b1}, introduced a new convolution in the special affine Fourier domain that works well with both the SAFT and its inverse leading to an analogue of the convolution and product formulas for the Fourier transform. They also introduced a second convolution that eliminates the phase factor in the convolution proposed by Xiang and Qin \cite{b2}.

\parindent=8mm \vspace{.1in}
Uncertainty principles are mathematical results that give limitations on the simul- taneous concentration of a function and its quaternion Fourier transform. They have implications in two main areas: quantum physics and signal analysis \cite{l1,l2,Fol,bk}. In quantum physics, they tell us that a particle's speed and position cannot both be measured with infinite precision. In signal analysis, they tell us that if we observe a signal only for a finite period of time, we will lose information about the frequencies the signal consists of. There are many ways to get the statement about concentration precise. This principle has been extended to different setups  by various researchers \cite{cl1,bk,Cow,gg,ss,dd,n1,Ste2,Wil}.

\parindent=8mm \vspace{.1in}
Motivated and inspired by the above work, we introduce a new special affine wavelet transform based on the novel convolution introduced in \cite{b2} and we call it {\it novel special affine wavelet transform} (NSAWT). We first establish some fundamental properties including Moyal's principle, inversion formula and the range theorem. Some Heisenberg type inequalities and Pitt's inequality are established for SAFT and  Heisenberg uncertainity principle is  also derived for NSAWT.

\parindent=8mm \vspace{.1in}
The rest of the paper is tailored as follows.In section 2, we introduce a notion of novel special affine wavelet transform (NSAWT) and establish a relationship between special affine Fourier transform (SAFT) and the proposed  NSAWT. Section 3 is dedicated to the key harmonic analysis results related to novel special affine wavelet transform (NSAWT). Some generalizations of Heisenberg type inequalities are established in section 4.

\parindent=0mm \vspace{.1in}
{\bf{2. Novel Special Affine Wavelet Transform}}

\parindent=0mm \vspace{.1in}
In this section, we first notion of  novel special affine wavelet transform wich is based on novel convolution.  Furthermore we establish  a relationship between special affine Fourier transform (SAFT) and the proposed novel special affine wavelet transform (NSAWT).

\parindent=8mm \vspace{.1in}
Firstly, we recall  the following definition of novel convolution  and the corresponding convolution theorem given by Bhandari and Zayed \cite{b1}.

\parindent=0mm \vspace{.1in}
{\bf {Definition 2.1.}(Chirp Modulation)} Let  $\bf \Lambda_S$ be the augumented SAFT matix. We  the modulation function $m_{\bf\Lambda_S}$ as follows
\begin{align*}
m_{\bf\Lambda_S}(x)=\exp\left\{\dfrac{iAx^2}{2B}\right\}.\tag{2.1}
\end{align*}

\parindent=0mm \vspace{.0in}
Then, for  a given $f\in L^2(\mathbb R)$, the chirp modulated functions associated with the  augumented  SAFT matrix $\bf \Lambda_S$  and inverse SAFT matrix  ${\bf \Lambda_S}^{inv}$are defined in the following manner:
\begin{align*}
\overset{\rightharpoonup} f(x)=m_{\bf\Lambda_S}(x)\,f(x)\quad ;\quad \overset{\leftharpoonup} f(x)=m_{\bf\Lambda_S}^*(x)\,f(x).\tag{2.2}
\end{align*}
and
\begin{align*}
\overset{\rightharpoonup} f(x)=m_{{\bf \Lambda_S}^{inv}}(x)\,f(x)\quad ;\quad \overset{\leftharpoonup} f(x)=m^*_{{\bf \Lambda_S}^{inv}}(x)\,f(x).\tag{2.3}
\end{align*}

\parindent=8mm \vspace{.1in}

\parindent=0mm \vspace{.1in}
{\bf Definition 2.2.(SAFT Convolution)} Let  $f,g\in L^2(\mathbb R)$ be two given functions .The special affine convolution $\ast_{\bf \Lambda_S}$ is defined as
\begin{align*}
h(t)=\Big(f\ast_{\bf\Lambda_S} g\Big)(t)=K_{B}\,m^*_{\bf \Lambda_S}(t)\left(\overset{\rightharpoonup} f(x)\ast\overset{\rightharpoonup} g(t)\right),\tag{2.4}
\end{align*}
where $\ast$ denote the usual convolution operator.

\parindent=0mm \vspace{.1in}
{\bf Lemma 2.3. (SAFT Convolution Theorem)} Let $f$ and $g$ be two functions such that $h(t) = (f \ast_{\bf \Lambda_S} g) (t)$ exists, then            
\begin{align*}
\widehat h_{\bf \Lambda_S}(\omega)= \Phi_{\bf\Lambda_S}(\omega)\widehat f_{\bf\Lambda_S}(\omega)\,\widehat g_{\bf\Lambda_S}(\omega),\tag{2.5}
\end{align*}
where
\begin{align*}
\Phi_{\Lambda_S}(\omega)=\exp\left\{\dfrac{i\omega(Dp-Bq)}{B}\right\}\exp\left\{-\dfrac{iD\omega^2}{2B}\right\} .\tag{2.6}
\end{align*}

\parindent=8mm \vspace{.1in}
 On the  basis of   SAFT convolution defined in Definition 2.1, we shall introduce the notion of novel special affine wavelet transform (NSAWT).

\parindent=0mm \vspace{.1in}

{\bf Definition 2.4.} For any finite energy signal $f\in L^2(\mathbb R)$, the continuous novel  special affine wavelet transform of $f$ with respect to the wavelet  $\psi\in L^2(\mathbb R)$ is defined by
$${\mathcal{A}}^{\psi}_{\bf\Lambda_S} [f](t,\zeta)=\int_{\mathbb R} f(x)\, {\psi^*}_{t,\zeta}^{\bf\Lambda_S}(x)\, dx,\quad t\in\mathbb R,\, \zeta\in\mathbb R^+,\eqno(2.7)$$
where $\psi_{t,\zeta}^{\bf\Lambda_S}(x)$ is given as follows
$$\psi_{t,\zeta}^{\bf\Lambda_S}(x)= \dfrac{K_{B}}{\sqrt {\zeta}}\, \psi\left(\dfrac{x-t}{\zeta} \right)\exp\left\{\dfrac{iAx(t-x)}{B}\right\}.\eqno(2.8)$$

\parindent=0mm \vspace{.1in}
It is worth noting that the  NSAWT  boils down to some existing integral transforms as well as gives birth to some new time-frequency transforms as mentioned below:

\parindent=0mm \vspace{.1in}
(i)~ For ${\bf\Lambda_S}=\left(A,B,C,D: 0,0\right)$, we get a  a novel linear cannonical wavelet transform.

\parindent=0mm \vspace{.1in}
(ii)~ For ${\bf\Lambda_S}=\left(\cos\theta,\sin\theta,-\sin\theta,\cos\theta: p,q\right)$, $\theta\neq n\pi$, we get a  a novel fractional wavelet transform defined  by
\begin{align*}
\Big({\mathcal{A}}^{\psi}_{\bf\Lambda_S} f\Big)(t,\zeta)=\dfrac{1}{\sqrt {2\pi \xi \left|\sin\theta\right|}}\int_{\mathbb R} f(x)\, \psi^*\left(\dfrac{x-t}{\zeta }\right)\exp\big\{ (ix^2-itx)\cot\theta\big\}\,dx.
\end{align*}

\parindent=0mm \vspace{.1in}

(iii)~For ${\bf\Lambda_S}=(1,B,0,1: p,q)$, $B\neq 0$, we  obtain a novel Fresnel-wavelet transform:
  \begin{align*}
  \Big({\mathcal{A}}^{\psi}_{\Lambda_S} f\Big)(t,\zeta)=\dfrac{K_{B}}{\sqrt {\zeta}}\int_{\mathbb R} f(x)\,\psi^*\left(\dfrac{x-t}{\zeta} \right)\exp\left\{\dfrac{(ix^2-itx)}{B}\right\}\,dx.
  \end{align*}

\parindent=8mm \vspace{.0in}
Now, we proceed to establish a fundamental relationship between the special affine Fourier transform given by (1.1) and the proposed novel special affine wavelet transform  defined in Definition 2.4..

\parindent=0mm \vspace{.1in}

{\bf Theorem 2.5.}   Let ${\mathcal{A}}^{\psi}_{\bf\Lambda_S} [f](t,\zeta)$ and $\mathcal F_{SAFT}\big[f\big](\omega)$  be the continuous novel  special affine wavelet transform  and the special affine Fourier transform of any finite energy signal $f\in L^2(\mathbb R)$. Then, we have
{\small
\begin{align*}
\mathcal F_{SAFT}\Big[\big({\mathcal{A}}^{\psi}_{\Lambda_S} f\big)(t,\zeta)\Big](\omega)=\sqrt{\zeta} \exp\left\{\dfrac{i}{2B} \big(2\omega\zeta(Dp-Bq)-D\omega^2\zeta^2\big)\right\}\widehat f_{\bf\Lambda_{S}}(\omega) \widehat \Psi_{\bf\Lambda_{S}}(\zeta\omega,\zeta),\tag{2.9}
\end{align*}
where
\begin{align*}
\Psi(x,\zeta)=\exp\left\{\dfrac{i}{2B}\big(Ax^2(\zeta^2-1)+2xp(\zeta-1)\big)\right\}\tilde \psi(x).\tag{2.10}
\end{align*}

\parindent=0mm \vspace{.0in}
{\it Proof.} By the definition of special affine Fourier transform, we have
\begin{align*}
\mathcal F_{SAFT}\Big[\big({\mathcal{A}}^{\psi}_{\Lambda_S} f\big)(t,\zeta)\Big](\omega)&= \mathcal F_{SAFT}\left[\left\{ f(x) \ast_{\bf\Lambda_S} \dfrac{1}{\sqrt{\zeta}}\psi^*\left(\frac{-x}{\zeta}\right) \exp\left\{\dfrac{-iAx^2}{2B}\right\}\right\} \right](\omega)\\\\
&=\exp\left\{\dfrac{i\omega(Dp-Bq)}{B}\right\}\exp\left\{-\dfrac{iD\omega^2}{2B}\right\}(\omega)\\\
&\qquad\qquad\times\mathcal F_{SAFT}\big[f(x)\big](\omega)\,\mathcal F_{SAFT}\Big[\psi^*\left(\dfrac{-x}{\zeta}\right)\Big](\omega),\tag{2.11}
\end{align*}
Further, we have
\begin{align*}
&\mathcal F_{SAFT}\left[\dfrac{1}{\sqrt{\zeta}}\psi^*\left(\frac{-x}{\zeta}\right) \exp\left\{\dfrac{-iAx^2}{2B}\right\}\right](\omega)\\\\
&\quad=\dfrac{K_B}{\sqrt{\zeta}}\int_{\mathbb R} \psi^*\left(-\dfrac{x}{\zeta}\right)\, \exp\left\{\dfrac{i}{2B}\big( Ax^2+D\omega^2+2x(p-\omega)-2\omega(Dp-Bq)\big)\right\}dx\\\\
&\quad=K_{B}\sqrt{\zeta}\int_{\mathbb R} \psi^*(-z)\,\exp\left\{\dfrac{i}{2B}\big( Az^2\zeta^2+D\omega^2-2\omega(Dp-Bq)\big)\right\}\exp\left\{\dfrac{i z\zeta(p-\omega)}{B}\right\}dz
\end{align*}
\begin{align*}
&\quad=K_{B}\sqrt{\zeta} \exp\left\{\dfrac{i}{2B}\big(D\omega^2-2\omega(Dp-Bq)-D(\zeta\omega)^2+2\zeta\omega(Dp-Bq)\big)\right\} \\\\
&\quad\quad\qquad\qquad\qquad\times \int_{\mathbb R} \psi^*(-z)\,\exp\left\{\dfrac{i}{2B}\big(A\zeta^2z^2+2z\zeta p-Az^2-2zp\big)\right\}\, \dfrac{{\mathcal K}_{\bf\Lambda_{S}}(z,\zeta\omega)}{K_B}\,dz\\\\
&\quad=\sqrt{\zeta} \exp\left\{\dfrac{i}{2B}\big(D\omega^2(1-\zeta^2)+2\omega(Dp-Bq)(\zeta-1)\big)\right\} \\\
&\qquad\qquad\qquad\quad\qquad\times \int_{\mathbb R} \psi  (-z)\,\exp\left\{\dfrac{i}{2B}\big(Az^2(\zeta^2-1)+2zp(\zeta-1)\big)\right\}\,{\mathcal K}_{\bf\Lambda_{S}}(z,\zeta\omega)\,dz\\\\
&\quad=\sqrt{\zeta} \exp\left\{\dfrac{i}{2B}\big(D\omega^2(1-\zeta^2)+2\omega(Dp-Bq)(\zeta-1)\big)\right\}\,\widehat\Psi_{\bf\Lambda_S}(\zeta\omega,\zeta),
\end{align*}
where
\begin{align*}
\Psi(x,\zeta)=\exp\left\{\dfrac{i}{2B}\big(Ax^2(\zeta^2-1)+2xp(\zeta-1)\big)\right\}  \psi^*(-x).\tag{2.12}
\end{align*}
From the equation  (2.11), we obtain the required result as
\begin{align*}
\mathcal F_{SAFT}\Big[\big({\mathcal{A}}^{\psi}_{\Lambda_S} f\big)(t,\zeta)\Big](\omega)&=\sqrt{\zeta} \exp\left\{\dfrac{i}{2B}\big(D\omega^2(1-\zeta^2)+2\omega(Dp-Bq)(\zeta-1)\big)\right\}\\\
&\qquad\qquad\times\exp\left\{ \dfrac{i\omega(Dp-Bq)}{B}-\dfrac{iD\omega^2}{2B}\right\}\widehat f_{\bf\Lambda_S}(\omega)\widehat \Psi_{\bf\Lambda_{S}}(\zeta\omega,\zeta)\\\\
&=\sqrt{\zeta} \exp\left\{\dfrac{i}{2B}\big(2\omega\zeta(Dp-Bq)-D\omega^2\zeta^2\big)\right\}\widehat f_{\Lambda_{S}}(\omega) \widehat \Psi_{\bf\Lambda_{S}}(\zeta\omega,\zeta).
\end{align*}
This completes the proof.\quad\fbox

\parindent=8mm \vspace{.1in}

From the Theorem 2.5, we conclude that if the analyzing functions $\psi_{t,\zeta}^{\bf\Lambda_S}(x)$ are supported in the time-domain or the special affine Fourier domain, then the proposed transform $\mathcal{W}^M_\psi \big[f \big](a,b)$ is accordingly supported in the respective domains. This implies that the special affine wavelet transform is capable of providing the simultaneous information of the time and the special affine frequency in the time-frequency domain. To be more specific, suppose that $\psi(t)$ is the window with centre $E_{\psi}$ and radius $\Delta_{\psi}$ in the time domain. Then, the centre and radii of the time-domain window function $\psi_{t,\zeta}^{\bf\Lambda_S}(x)$  of the proposed transform (2.7) is given by
\begin{align*}
E\Big[\psi_{t,\zeta}^{\bf\Lambda_S}(x)\Big]&=\dfrac{\displaystyle\int_{-\infty}^{\infty}x\Big|\psi_{t,\zeta}^{\bf\Lambda_S}(x)\Big|^2 dx}{\displaystyle\int_{-\infty}^{\infty}\Big|\psi_{t,\zeta}^{\bf\Lambda_S}(x)\Big|^2 dx}\\\\
&=\dfrac{\displaystyle\int_{-\infty}^{\infty}x\Big|\psi_{t,\zeta}(x)\Big|^2 dx}{\displaystyle\int_{-\infty}^{\infty}\Big|\psi_{t,\zeta}(x)\Big|^2dx}\\\\
&=E\Big[\psi_{t,\zeta}(x)\Big]=t+\zeta E_{\psi}
\end{align*}
and
\begin{align*}
\Delta\Big[\psi_{t,\zeta}^{\bf\Lambda_S}(x)\Big]&=\left\{\dfrac{\displaystyle\int_{-\infty}^{\infty}\Big(x-(t+\zeta E_{\psi})\Big)\Big|\psi_{t,\zeta}^{\bf\Lambda_S}(x)\Big|^2 dx}{\displaystyle\int_{-\infty}^{\infty}\left|\psi_{t,\zeta}^{\bf\Lambda_S}(x)\right|^2 dx}\right\}^{1/2}\\\\
&=\left\{\dfrac{\displaystyle\int_{-\infty}^{\infty}\Big(x-t-\zeta E_{\psi}\Big)\Big|\psi_{t,\zeta}(x)\Big|^2 dx}{\displaystyle\int_{-\infty}^{\infty}\Big|\psi_{t,\zeta}(x)\Big|^2dt}\right\}^{1/2}\\\\
&=\Delta\Big[\psi_{t,\zeta}(x)\Big]= \zeta\Delta_{\psi},
\end{align*}
respectively. Let $\Gamma (\omega)$ be the window function in the special affine Fourier  domain (SAFD) given by
\begin{align*}
\Gamma(\omega)=\mathcal{F}_{SAFT}\left[\exp\left\{\dfrac{i}{2B}\Big(2x\zeta p-2xp-Ax^2\Big) \right\}\psi^*(-x)\right]\left(\omega\right).
\end{align*}

\parindent=0mm \vspace{.0in}
Then, we can derive the center and radius of the special affine Fourier domain (SAFD) window function
\begin{align*}
\Gamma(\zeta \omega)=\mathcal{F}_{SAFT}\left[\exp\left\{\dfrac{i}{2B}\Big(2x\zeta p-2xp-Ax^2)\Big) \right\}\psi^*(-x)\right]\left(\zeta \omega\right)
\end{align*}
appearing in (2.12) as
\begin{align*}
E\big[\Gamma \left(\zeta \omega\right)\big]= \dfrac{\displaystyle\int_{-\infty}^{\infty}(a\omega)\big|\Gamma(\zeta \omega)\big|^2 d\omega}{\displaystyle\int_{-\infty}^{\infty}\big|\Gamma (\zeta \omega)\big|^2 d\omega} = \zeta \, E_\Gamma,
\end{align*}
and
\begin{align*}
\Delta\big[\Gamma\left(\zeta \omega \right)\big] =\zeta \,\Delta_\Gamma. 
\end{align*}
Thus, the $Q$-factor of the  proposed transform (2.7) is given by
\begin{align*}
Q_{NSAWT}=\dfrac{\text{width of the window function }}{\text{centre of the window function }}=\dfrac{\Delta\big[\Gamma\left(\zeta \omega\right)\big]}{E\big[\Gamma\left(\zeta \omega \right)\big]}=\dfrac{\Delta_\Gamma}{E_\Gamma}=\text{constant},
\end{align*}
which is independent of the uni-modular matrix  ${\bf\Lambda_S}$ and the scaling parameter $\zeta$.  Therefore, the localized time and frequency characteristics of the novel special affine wavelet transform (NSAWT)  are given in the  time and frequency windows
\begin{align*}
\Big[t+\zeta E_{\psi}-\zeta \Delta_{\psi},~t+\zeta E_{\psi}+\zeta \Delta_{\psi}\Big] \quad \text{and}\quad \Big[a\,E_\Gamma-\zeta\,\Delta_\Gamma,~\zeta\,E_\Gamma+\zeta\,\Delta_\Gamma\Big],
\end{align*}
respectively. Hence, the joint resolution of the continuous  novel special affine wavelet transform (NSAWT)  in the time-frequency domain is described by a flexible window $\psi$ having a total spread  $4\Delta_{\psi}\Delta_\Gamma$ and is given by
\begin{align*}
\Big[t+\zeta E_{\psi}-\zeta \Delta_{\psi},~t+\zeta E_{\psi}+\zeta \Delta_{\psi}\Big]  \times  \Big[a\,E_\Gamma-\zeta\,\Delta_\Gamma,~\zeta\,E_\Gamma+\zeta\,\Delta_\Gamma\Big],
\end{align*}

\parindent=0mm \vspace{.1in}
{\bf{3.  Basic Properties of  Novel Special Affine Wavelet Transform}}

\parindent=0mm \vspace{.1in}
In this section, we establish  fundamental properties of the novel special affine wavelet transform (NSAWT). Some well known harmonic analyis results namely  Moyal's principle, inversion formula,  characterization of the range of the novel special aafine wavelet transform are derived.

\parindent=8mm \vspace{.1in}
Now, we proceed to state  some fundamental properties of the novel special affine wavelet transform (NSAWT)  defined in Definition 2.4.

\parindent=0mm \vspace{.1in}

{\bf Theorem 3.1.}  For any functions  $f,g\in L^2(\mathbb R)$  and $\alpha,\beta, \gamma\in\mathbb R,~ \mu\in\mathbb R^+$,  the  continuous  novel special affine wavelet transform  satisfies the following properties:

\parindent=0mm \vspace{.1in}
(i)~~${\text {\bf Linearity:}}~~\qquad{\mathcal{A}}^{\psi}_{\Lambda_S} \big[\alpha f+\beta g\big] (t,\zeta)= \alpha~{\mathcal{A}}^{\psi}_{\bf\Lambda_S} \big[f\big](t,\zeta)+\beta~{\mathcal{A}}^{\psi}_{\bf\Lambda_S}\big[g\big](t,\zeta)$.

\parindent=0mm \vspace{.1in}
(ii)~${\text {\bf Translation:}}\,~ {\mathcal{A}}^{\psi}_{\bf\Lambda_S}\big[f(x-\gamma)\big](t,\zeta)=\exp\left\{-\dfrac{iA\gamma(t-\gamma)}{B}\right\}{\mathcal{A}}^{\psi}_{\bf\Lambda_S}\left[\exp\left\{\dfrac{iAy\gamma}{B}\right\} f(y)\right](t-\gamma, \zeta).$

\parindent=0mm \vspace{.1in}
(iii)~${\text {\bf Scaling:}}~\quad{\mathcal{A}}^{\psi}_{\bf\Lambda_S}\big[f(\gamma x)\big](t,\zeta)=\dfrac{1}{\sqrt{\gamma}}{\mathcal{A}}^{\psi}_{{\bf\Lambda^{\prime}_S}}\big[f\big](\gamma t,\gamma \zeta),$
where
$${\bf\Lambda_S^{\prime}}=\left[\begin{array}{ccc}\phantom{-}A &~ \gamma^2 B :& p\\ ~C & \phantom{-}D  : &q \\ \end{array}\right]$$

\parindent=0mm \vspace{.1in}

{\it Proof.} These properties are obvious, therefore we omit the proofs.

\parindent=8mm \vspace{.1in}

Now, we shall study some important theorems including the Moyal's theorem, inversion formula and range theorem pertaining to the novel special affine wavelet transform defined in Definition 2.4.  Firstly, we shall derive the admissibility condition associated with the novel special affine wavelet transform.

\parindent=0mm \vspace{.1in}
{\bf Theorem 3.2 {(\bf  Admissibility Condition).}} { Let $\psi\in L^2(\mathbb R)$  be a given function, then $\psi$ is said to be admissible if
\begin{align*}
C_{\psi}=\int_{\mathbb R^+} \dfrac{\Big|\widehat \Psi_{\bf\Lambda_{S}}(\zeta\omega,\zeta)\Big|^2}{\zeta} \, d\zeta<\infty,\quad a.e.\quad \omega \in\mathbb R\tag{3.1}
\end{align*}
where $\Psi$ is given by (2.10).}

\parindent=0mm \vspace{.1in}
{\it Proof.} For any  $f \in L^2(\mathbb R)$, we have
\begin{align*}
\int_{\mathbb R\times \mathbb R^+}\Big|\Big\langle f, \psi_{t,\zeta}^{\bf\Lambda_S}\Big\rangle\Big|^2 \dfrac{dt\,d\zeta}{\zeta^2} &= \int_{\mathbb R\times \mathbb R^+}\Big|\Big( f(x) \ast_{\bf\Lambda_S} \dfrac{1}{\sqrt{\zeta}}\psi^*\left(\frac{-x}{\zeta}\right) \exp\left\{\dfrac{-iAx^2}{2B}\right\} \Big)(t)\Big|^2 \dfrac{dt\,d\zeta}{\zeta^2}\\\\
&= \int_{\mathbb R\times \mathbb R^+}\Big|\mathcal F_{SAFT}\left[\Big( f(x) \ast_{\bf\Lambda_S} \dfrac{1}{\sqrt{\zeta}}\psi^*\left(\frac{-x}{\zeta}\right) \exp\left\{\dfrac{-iAx^2}{2B}\right\} \Big)\right](\omega)\Big|^2 \dfrac{d\omega\,d\zeta}{\zeta^2}\\\\
&= \int_{\mathbb R\times \mathbb R^+}\big|\zeta\big|\left|\widehat f_{\bf\Lambda_{S}}(\omega)\right|^2\Big|\widehat \Psi_{\bf\Lambda_{S}}(\zeta\omega,\zeta)\Big|^2 \dfrac{d\omega\,d\zeta}{\zeta^2}\\\\
&= \int_{\mathbb R}\left|\widehat f_{\bf\Lambda_{S}}(\omega)\right|^2 \left\{ \int_{\mathbb R^+} \dfrac{\Big|\widehat \Psi_{\bf\Lambda_{S}}(\zeta\omega,\zeta)\Big|^2}{\zeta}\,d\zeta\right\} d\omega.\tag{3.2}
\end{align*}
On putting $f=\psi$,~ (3.2) reduces  to
\begin{align*}
\int_{\mathbb R\times \mathbb R^+}\Big|\Big\langle \psi, \psi_{t,\zeta}^{\bf\Lambda_S}\Big\rangle\Big|^2 \dfrac{dt\,d\zeta}{\zeta^2} =  \int_{\mathbb R}\left|\widehat \psi_{\bf\Lambda_{S}}(\omega)\right|^2 \left\{ \int_{\mathbb R^+} \dfrac{\Big|\widehat \Psi_{\bf\Lambda_{S}}(\zeta\omega,\zeta)\Big|^2}{\zeta}\,d\zeta\right\} d\omega.\tag{3.3}
\end{align*}
Since $\psi\in L^2(\mathbb R)$, therefore we conclude that the R.H.S of (3.3) is finite provided
\begin{align*}
C_{\psi}=\int_{\mathbb R^+} \dfrac{\Big|\widehat \Psi_{\bf\Lambda_{S}}(\zeta\omega,\zeta)\Big|^2}{\zeta} \, d\zeta<\infty,\quad a.e. \quad \omega\in\mathbb R.\quad\square
\end{align*}

\parindent=8mm \vspace{.1in}

The following is the Moyal's principle for the novel special affine wavelet transform (NSAWT).

\parindent=0mm \vspace{.1in}
{\bf Theorem 3.3 {(\bf Moyal's Principle).}}   Let ${\mathcal{A}}^{\psi}_{\bf\Lambda_S}\big[f\big](t,\zeta)$  and ${\mathcal{A}}^{\psi}_{\bf\Lambda_S}\big[g\big](t,\zeta)$  be the  novel special affine wavelet transforms of  $f$ and $g$ belonging to $L^2(\mathbb R)$, respectively. Then, we have
\begin{align*}
\int_{\mathbb R\times \mathbb R^+}{\mathcal{A}}^{\psi}_{\bf\Lambda_S}\big[f\big](t,\zeta)\,{\mathcal{A}}^{\star \psi }_{\bf\Lambda_S}\big[g\big](t,\zeta)\, \dfrac{dt\,d\zeta}{\zeta^2}= C_{\psi}\,\Big\langle f,\,g \Big\rangle_{L^2(\mathbb R)},\tag{3.4}
\end{align*}
where $C_{\psi}$ is given by (3.1).

 \parindent=0mm \vspace{.1in}
{\it Proof.} Applying  Theorem 2.5, we have for any pair of square integrable functions $f$ and $g$
\begin{align*}
{\mathcal{A}}^{\psi}_{\bf\Lambda_S}\big[f\big](t,\zeta)&=\sqrt{\zeta}\int_{\mathbb R} \exp\left\{\dfrac{i}{2B}\big(2\omega\zeta(Dp-Bq)-D\omega^2\zeta^2\big)\right\}\widehat f_{\bf\Lambda_{S}}(\omega) \widehat \Psi_{\bf\Lambda_{S}}(\zeta\omega,\zeta)\, \mathcal{K}_{\bf\Lambda_S^{inv}}(\omega,t)\, d\omega
\end{align*}
and
\begin{align*}
{\mathcal{A}}^{\psi}_{\bf\Lambda_S}\big[g\big](t,\zeta)&=\sqrt{\zeta}\int_{\mathbb R} \exp\left\{\dfrac{i}{2B}\big(2\eta\zeta(Dp-Bq)-D\eta^2\zeta^2\big)\right\}\widehat g_{\bf\Lambda_{S}}(\eta) \widehat \Psi_{\bf\Lambda_{S}}(\zeta\eta,\zeta)\, \mathcal{K}_{\bf\Lambda_S^{inv}}(\eta,t)\, d\eta,
\end{align*}
where $\Psi$ is given by (2.10), respectively. Further, we have
\begin{align*}
&\int_{\mathbb R\times \mathbb R^+}{\mathcal{A}}^{\psi}_{\bf\Lambda_S}\big[f\big](t,\zeta)\,{\mathcal{A}}^{\star\psi}_{\bf\Lambda_S}\big[g\big](t,\zeta)\, \dfrac{dt\,d\zeta}{\zeta^2}\\\\
&\qquad=\int_{\mathbb R\times \mathbb R\times \mathbb R\times \mathbb R^+}\exp\left\{\dfrac{i}{2B}\big(2\zeta(Dp-Bq)(\omega-\eta)-D\zeta^2(\omega^2-\eta^2)\big) \right\} \\\\
&\qquad\qquad\qquad\times~\widehat f_{\bf\Lambda_S}(\omega)\,\widehat g^*_{\bf\Lambda_S}(\eta)\, \widehat \Psi_{\bf\Lambda_{S}}(\zeta\omega,\zeta)\,\widehat \Psi^\star_{\bf\Lambda_{S}}(\zeta\eta,\zeta)\,\mathcal{K}_{\bf\Lambda_S^{inv}}(\omega,t)\,\mathcal{K}^\star_{\bf\Lambda_S^{inv}}(\eta,t)\, \dfrac{dt\,d\omega\,d\eta\,d\zeta}{\zeta}\\\\
\end{align*}
\begin{align*}
&\qquad= \int_{\mathbb R\times \mathbb R\times \mathbb R^+}\exp\left\{\dfrac{i}{2B}\big(2\zeta(Dp-Bq)(\omega-\eta)-D\zeta^2(\omega^2-\eta^2)\big) \right\} \\\\
&\qquad\qquad\qquad\times~\widehat f_{\Lambda_S}(\omega)\,\widehat g^\star_{\bf\Lambda_S}(\eta) \widehat \Psi_{\bf\Lambda_{S}}(\zeta\omega,\zeta)\,\widehat \Psi^\star_{\bf\Lambda_{S}}(\zeta\eta,\zeta)\left\{\int_{\mathbb R}\mathcal{K}_{\bf\Lambda_S^{inv}}(\omega,t)\,\mathcal{K}^\star_{\bf\Lambda_S^{inv}}(\eta,t)\,dt\right\} \dfrac{d\omega\,d\eta\,d\zeta}{\zeta}\\\\
&\qquad=\int_{\mathbb R\times \mathbb R\times \mathbb R^+}\exp\left\{\dfrac{i}{2B}\big(2\zeta(Dp-Bq)(\omega-\eta)-D\zeta^2(\omega^2-\eta^2)\big) \right\} \\\\
&\qquad\qquad\qquad\qquad\qquad\times~\widehat f_{\Lambda_S}(\omega)\,\widehat g^\star_{\bf\Lambda_S}(\eta)  \widehat \Psi_{\Lambda_{S}}(\zeta\omega,\zeta)\,\widehat \Psi^\star_{\Lambda_{S}}(\zeta\eta,\zeta)\,\delta(\omega-\eta)\,\dfrac{d\omega\,d\eta\,d\zeta}{\zeta}\\\\
&\qquad=\int_{\mathbb R\times \mathbb R^+ } \widehat f_{\bf\Lambda_S}(\omega)\,\widehat g^\star_{\bf\Lambda_S} (\omega)\,\Big|\widehat \Psi_{\bf\Lambda_{S}}(\zeta\omega,\zeta)\Big|^2\,\dfrac{d\omega\,d\zeta}{\zeta}\\\\
&\qquad=\int_{\mathbb R} \widehat f_{\bf\Lambda_S}(\omega)\,\widehat g^\star_{\bf\Lambda_S} (\omega) \left\{ \int_{\mathbb R^+} \dfrac{\Big|\widehat \Psi_{\bf\Lambda_{S}}(\zeta\omega,\zeta)\Big|^2}{\zeta} \, d\zeta\right\}d\omega\\\
&\qquad=C_{\psi}\,\left\langle \widehat f_{\bf\Lambda_S},\, \widehat g_{\bf\Lambda_S}\right\rangle_{L^2(\mathbb R)}\\\
&\qquad= C_{\psi}\, \Big\langle f,\,g \Big\rangle_{L^2(\mathbb R)}.
\end{align*}
This completes the proof. \quad\fbox

\parindent=0mm \vspace{.1in}

It is worth to mention that for $f=g$, the above theorem reduces to :
\begin{align*}
\int_{\mathbb R\times \mathbb R^+}\left|{\mathcal{A}}^{\psi}_{\bf\Lambda_S}\big[f\big](t,\zeta)\right|^2 \dfrac{dt\,d\zeta}{\zeta^2} =C_{\psi}\,\big\|f\big\|^2_2.
\end{align*}
This is energy preserving theorem for the novel special affine wavelet transform.

\parindent=8mm \vspace{.1in}
Now, the following theorem is the inversion formula for the novel  special affine wavelet transform ${\mathcal{A}}^{\psi}_{\Lambda_S}\big[f\big](t,\zeta)$.

\parindent=0mm \vspace{.1in}
{\bf Theorem 3.4 (Inversion Formula).}  Let $f\in L^2(\mathbb R)$ be a given function and $\psi$ is admissible. If $~{\mathcal{A}}^{\psi}_{\Lambda_S}\big[f\big](t,\zeta)$ is the  novel  special affine wavelet transform of $f$, then $f$ can be reconstructed as
\begin{align*}
f(x)= \dfrac{1}{C_{\psi}} \int_{\mathbb R\times \mathbb R^+}{\mathcal{A}}^{\psi}_{\bf\Lambda_S}\big[f\big](t,\zeta)\,  \psi_{t,\zeta}^{\bf\Lambda_S}(x)\,\dfrac{dt\,d\zeta}{\zeta^2}, \quad a.e. \tag{3.5}
\end{align*}

\parindent=0mm \vspace{.0in}
{\it Proof.} By vitue of Moyal's principle, we have
\begin{align*}
\Big\langle f,\,g \Big\rangle&=\dfrac{1}{C_{\psi}}\int_{\mathbb R\times \mathbb R^+}{\mathcal{A}}^{\psi}_{\bf\Lambda_S}\big[f\big](t,\zeta)\,{\mathcal{A}^\star}^{\psi}_{\bf\Lambda_S}\big[g\big](t,\zeta)\, \dfrac{dt\,d\zeta}{\zeta^2}\\\\
&=\dfrac{1}{C_{\psi}}\int_{\mathbb R\times \mathbb R^+}{\mathcal{A}}^{\psi}_{\bf\Lambda_S}\big[f\big](t,\zeta)\,\left\{ \int_{\mathbb R} g^*(x)\,\psi_{t,\zeta}^{\bf\Lambda_S}(x) \,dx\right\} \dfrac{dt\,d\zeta}{\zeta^2}\\\\
&=\dfrac{1}{C_{\psi}}\int_{\mathbb R\times \mathbb R\times \mathbb R^+} {\mathcal{A}}^{\psi}_{\bf\Lambda_S}\big[f\big](t,\zeta)\,\psi_{t,\zeta}^{\bf\Lambda_S}(x)\, g^*(x)\, \dfrac{dt\,d\zeta}{\zeta^2}\\\\
&=\dfrac{1}{C_{\psi}}\left\langle \int_{\mathbb R\times \mathbb R^+} {\mathcal{A}}^{\psi}_{\bf\Lambda_S}\big[f\big](t,\zeta)\,\psi_{t,\zeta}^{\bf\Lambda_S}(x)\, \dfrac{dt\,d\zeta}{\zeta^2},\, g(x) \right\rangle.
\end{align*}
Since $g$ is chosen arbitrarily from $L^2(\mathbb R)$, therefore we obtain
\begin{align*}
f(x)= \dfrac{1}{C_{\psi}} \int_{\mathbb R\times \mathbb R^+} {\mathcal{A}}^{\psi}_{\bf\Lambda_S}\big[f\big](t,\zeta)\,\psi_{t,\zeta}^{\bf\Lambda_S}(x)\, \dfrac{dt\,d\zeta}{\zeta^2} \quad a.e.
\end{align*}

Thus the proof is completed.

\parindent=8mm \vspace{.1in}
The following theorem provides a complete characterization of the range of the novel special affine wavelet transform  ${\mathcal{A}}^{\psi}_{\Lambda_S}$.

\parindent=0mm \vspace{.1in}

{\bf Theorem 3.5} ({\bf Characterization of Range of ${\mathcal{A}}^{\psi}_{\Lambda_S}$}).  If $f\in L^2(\mathbb R\times \mathbb R^{+})$ and  $\psi$ is admissible wavelet, then $f$ belongs to the range  of  $ {\mathcal{A}}^{\psi}_{\bf\Lambda_S}$ if and only if 
\begin{align*}
f(t^{\prime},\zeta^{\prime})=\dfrac{1}{C_{\psi}}\int_{\mathbb R\times \mathbb R^{+}}f(t,\zeta)\,\Big\langle \psi_{t,\zeta}^{\bf\Lambda_S}, \psi_{t^{\prime},\zeta^{\prime}}^{\bf\Lambda_S}\Big\rangle_2\,\dfrac{dt\,d\zeta}{\zeta^2}.\tag{3.6}
\end{align*}
\parindent=0mm \vspace{.1in}
{\it Proof.} Let  $f$ belongs to  range of $ {\mathcal{A}}^{\psi}_{\Lambda_S}$. Then, there exists a square integrable function $g$, such that ${\mathcal{A}}^{\psi}_{\Lambda_S}g=f$. In order to show that $f$ satisfies (3.6), we proceed as
\begin{align*}
f\big(t^{\prime},\zeta^{\prime})&={\mathcal{A}}^{\psi}_{\bf\Lambda_S}\big[g\big](t^{\prime},\zeta^{\prime})\\\\
&=\int_{\mathbb R}g(x)\,{\psi^\star}_{t^{\prime},\zeta^{\prime}}^{\bf\Lambda_S}(x)\,dx\\\\
&=\dfrac{1}{C_{\psi}}\int_{\mathbb R}\left\{\int_{\mathbb R\times \mathbb R^+} {\mathcal{A}}^{\psi}_{\bf\Lambda_S}\big[g\big](t,\zeta)\,\psi_{t,\zeta}^{\bf\Lambda_S}(x)\, \dfrac{dt\,d\zeta}{\zeta^2}\right\}\,{\psi^\star}_{t^{\prime},\zeta^{\prime}}^{\bf\Lambda_S}(x)\,dx\\\\
&=\dfrac{1}{C_{\psi}}\int_{\mathbb R\times \mathbb R^+}{\mathcal{A}}^{\psi}_{\bf\Lambda_S}\big[g\big](t,\zeta)\,\left\{\int_{\mathbb R}\psi_{t,\zeta}^{\bf\Lambda_S}(x)\,
{\psi^\star}_{t^{\prime},\zeta^{\prime}}^{\bf\Lambda_S}(x)\,dx\right\}\,\dfrac{dt\,d\zeta}{\zeta^2}\\\\
&=\dfrac{1}{C_{\psi}}\int_{\mathbb R\times \mathbb R^{+}}f(t,\zeta)\,\Big\langle \psi_{t,\zeta}^{\bf\Lambda_S}, \psi_{t^{\prime},\zeta^{\prime}}^{\bf\Lambda_S}\Big\rangle_2\,\dfrac{dt\,d\zeta}{\zeta^2},
\end{align*}
which proves  the necessary part. Conversely, suppose that a square integrable function $f$ satisfies (3.6). In order to prove that $f$ belongs to range of $ {\mathcal{A}}^{\psi}_{\bf\Lambda_S}$, we need a function $g\in L^2(\mathbb R)$ satisfying ${\mathcal{A}}^{\psi}_{\bf\Lambda_S}g=f$. This required  function $g$ is constructed as follows
\begin{align*}
g(x)=\dfrac{1}{C_{\psi}}\int_{\mathbb R\times \mathbb R^+} f(t,\zeta)\,\psi_{t,\zeta}^{\bf\Lambda_S}(x)\, \dfrac{dt\,d\zeta}{\zeta^2}.\tag{3.7}
\end{align*}
It is clear that $\big\|g\big\|_2\leq\big\|f\big\|_2<\infty$; that is $g\in L^2(\mathbb R)$. Also, we have
\begin{align*}
{\mathcal{A}}^{\psi}_{\bf\Lambda_S}\big[g\big](t^{\prime},\zeta^{\prime})&=\int_{\mathbb R}g(x)\,{\psi^\star}_{t^{\prime},\zeta^{\prime}}^{\bf\Lambda_S}(x)\,dx\\\\
&=\dfrac{1}{C_{\psi}}\int_{\mathbb R}\left\{\int_{\mathbb R\times \mathbb R^+} f(t,\zeta)\,\psi_{t,\zeta}^{\bf\Lambda_S}(x)\, \dfrac{dt\,d\zeta}{\zeta^2}\right\}\,
{\psi^\star}_{t^{\prime},\zeta^{\prime}}^{\bf\Lambda_S}(x)\,dx\\\\
&=\dfrac{1}{C_{\psi}}\int_{\mathbb R\times \mathbb R^{+}}f(t,\zeta)\,\Big\langle \psi_{t,\zeta}^{\bf\Lambda_S}, \psi_{t^{\prime},\zeta^{\prime}}^{\bf\Lambda_S}\Big\rangle_2\,\dfrac{dt\,d\zeta}{\zeta^2}\\\\
&=f(t^{\prime},\zeta^{\prime}).
\end{align*}

This completes the proof. \qquad \fbox

\parindent=0mm \vspace{.1in}

{\bf Corollary 3.6} {(\bf Reproducing Kernel Hilbert Space).} For any admissible wavelet $\psi\in L^2(\mathbb R)$, the range of  ${\mathcal{A}}^{\psi}_{\Lambda_S}$ is a reproducing kernel Hilbert space embedded as a subspace in $L^2\big(\mathbb R\times\mathbb R^+\big)$ with the kernel given by
\begin{align*}
K_{\bf\Lambda_S}^{\psi}\big(t,\zeta;t^{\prime},\zeta^{\prime}\big)=\Big\langle \psi_{t,\zeta}^{\bf\Lambda_S}, \psi_{t^{\prime},\zeta^{\prime}}^{\bf\Lambda_S}\Big\rangle.\tag{3.8}\\
\end{align*}

\newpage
\parindent=0mm \vspace{.1in}
{\bf{ 4. Uncertainty Principles Associated with SAFT and NSAWT}}

\parindent=0mm \vspace{.1in}
The Pitt's inequality in the Fourier domain expresses a fundamental relationship between a sufficiently smooth function and the corresponding Fourier transform \cite{bk}. For every $f\in \mathbb S(\mathbb R)\subseteq L^2(\mathbb R)$, the inequality states that
\begin{align*}
\int_{\mathbb R} |\omega|^{-\alpha}\left|\mathcal F\big[f\big](\omega)\right|^2d{\omega}\le C_{\alpha}\int_{\mathbb R}\left|x\right|^{\alpha}\big|f(x)\big|^2dx,\quad 0\le\alpha<1
\end{align*}
where
\begin{align*}
C_{\alpha}=\pi^{\alpha}\left[\Gamma\left(\frac{1-\alpha}{4}\right)/\, \Gamma\left(\frac{1+\alpha}{4}\right)\right]^2,
\end{align*}
and $\Gamma(\cdot)$ denotes the well known Euler's gamma function.The Schwartz class in $L^2(\mathbb R)$ is defined by
\begin{align*}
\mathbb S\left(\mathbb R\right)=\left\{f\in C^{\infty}(\mathbb R): \sup_{t\in\mathbb R}t^{\beta}\,{\mathbb D}^{\gamma}f(t)<\infty\right\},
\end{align*}
where  $C^{\infty}(\mathbb R)$ is the class of smooth functions, $\beta,\gamma$ are non-negative integers, and ${\mathbb D}$ denotes the usual  differential operator.

\parindent=8mm \vspace{.1in}
For any  $f\in L^2(\mathbb R)$, the Heisenberg's uncertainty inequality in the special affine Fourier domain is given by \cite{Ste2}
\begin{align*}
\left\{\int_{\mathbb R} x^2\big|f(x)\big|^2 dx \right\}^{1/2}\left\{\int_{\mathbb R} {\omega}^2\left|\widehat f_{\bf\Lambda_{S}}(\omega)\right|^2 d\omega \right\}^{1/2}\ge \dfrac{|B|}{2}\left\{\int_{\mathbb R} \big|f(x)\big|^2 dx\right\},\tag{4.1}
\end{align*}
with equality if and only if $f$ is a multiple of a suitable Gaussian function. 

\parindent=8mm \vspace{.1in}
 By adopting the strategy analogous to Wilcok \cite{Wil} and Cowling and Price \cite{Cow}, we establish a generalization of  the uncertainity principle given by (4.1). Furthermore, we  derive an uncertainty inequality comparing the localization of the special affine Fourier transform (SAFT) of a function $f$ with the novel special affine wavelet transform (NSAWT)  ${\mathcal{A}}^{\psi}_{\Lambda_S}\big[f\big](t,\zeta)$, regarded as a function of the time variable $t$.

\parindent=0mm \vspace{.1in}
{\bf Theorem 4.1.} {For any $f\in L^2(\mathbb R)$, the generalized uncertainty inequality for the special affine Fourier transform (1.1) is given by:
\begin{align*}
\left\{\int_{\mathbb R} \big|x\big|^p\big|f(x)\big|^p dx \right\}^{1/p}\left\{\int_{\mathbb R} \big|{\omega}\big|^p\left|\widehat f_{\bf\Lambda_{S}}(\omega)\right|^p d\omega \right\}^{1/p}\ge \dfrac{|B|^{(p+2)/2p}}{2}\Big\|f\Big\|^2_2,\quad 1\le p\le 2.\tag{4.2}
\end{align*}

\parindent=0mm \vspace{.1in}
{\it Proof.} For any $f\in L^2(\mathbb R)$, the generalized uncertainty inequality in the classical Fourier domain is given by
\begin{align*}
\left\{\int_{\mathbb R} \big|x\big|^p\big|f(x)\big|^p dx \right\}^{1/p}\left\{\int_{\mathbb R} \big|{\omega}\big|^p\left|\mathcal F\big[f\big](\omega)\right|^p d\omega \right\}^{1/p}\ge \dfrac{1}{2}\left\{\int_{\mathbb R} \big|f(x)\big|^2 dx\right\},~ 1\le p\le 2\tag{4.3}
\end{align*}
where $\mathcal F\big[f\big]$ denotes the classical Fourier transform of $f$.
We rewrite the definition of the special affine Fourier transform (1.1) as
\begin{align*}
\widehat f_{\bf\Lambda_{S}}(\omega)=K_B\int_{\mathbb R}g(x)\,\exp\left\{\dfrac{i}{2B}\left(D\omega^2-2x\omega-2\omega(Dp-Bq)\right)\right\}\,dx,\tag{4.4}
\end{align*}
where $g(x)=f(x)\exp\big\{i\big(Ax^2+2xp\big)/2B\big\}$. From equation (4.4), it is quite evident that
\begin{align*}
\widehat f_{\bf\Lambda_{S}}(\omega)
=\dfrac{1}{\sqrt{|B|}}\exp\left\{\dfrac{i}{2B}\left(D\omega^2-2\omega(Dp-Bq)\right)\right\}\,\mathcal F\Big[g\Big]\left(\dfrac{\omega}{B}\right),\tag{4.5}
\end{align*}
so that, $\sqrt{|B|}\,\big|\widehat f_{\bf\Lambda_{S}}(B\omega)\big|=\big|\mathcal F\big[g\big]\left(\omega\right)\big|$. Invoking (4.3) for the function $g$, we obtain the generalized uncertainty inequality for the special affine Fourier transform:
\begin{align*}
\left\{\int_{\mathbb R} \big|x\big|^p\big|f(x)\big|^p dx \right\}^{1/p}\left\{\int_{\mathbb R} \big|{\omega}\big|^p\left|\widehat f_{\bf\Lambda_{S}}(\omega)\right|^p d\omega \right\}^{1/p}\ge \dfrac{|B|^{(p+2)/2p}}{2}\Big\|f\Big\|^2_2,\quad 1\le p\le 2.
\end{align*}

\parindent=0mm \vspace{.1in}

{\it Remark:} For $p=2$, the generalized uncertainty principle (4.2) boils down to the classical Heisenberg's uncertainty principle for the special affine Fourier transform.

\parindent=8mm \vspace{.1in}
In view of   classical Pitt's inequality and the relationship between the classical Fourier and the special affine Fourier transform  given by (4.5), one clears obtain the following Pitt's inequality for SAFT.

\parindent=0mm \vspace{.1in}
{\bf Theorem 4.2 (Pitt's Inequality).}  For any $f\in \mathbb S(\mathbb R)$, the Pitt's inequality for the  special affine Fourier transform (1.1) is given by:
\begin{align*}
|B|^{\alpha}\int_{\mathbb R} |\omega|^{-\alpha}\big|\widehat f_{\bf\Lambda_{S}}(\omega)\big|^2d{\omega}\le C_{\alpha}\int_{\mathbb R}\left|x\right|^{\alpha}\big|f(x)\big|^2dx,\quad 0\le\alpha<1.
\end{align*}

\parindent=8mm \vspace{.1in}
Now, we shall derive an uncertainty inequality governing the simultaneous localization of ${\mathcal F}_{SAFT}\big[f\big](\omega)$ and ${\mathcal{A}}^{\psi}_{\Lambda_S}\big[f\big](t,\cdot)$.

\parindent=0mm \vspace{.1in}
{\bf Theorem 4.3.}  If ${\mathcal{A}}^{\psi}_{\bf\Lambda_S}\big[f\big](t,\zeta)$ is the novel special affine wavelet transform of any nontrivial function  $f\in L^2(\mathbb R)$, then the following uncertainty inequality holds:
{\small
\begin{align*}
&\left\{\int_{\mathbb R\times\mathbb R^+} t^2\Big| {\mathcal{A}}^{\psi}_{\bf\Lambda_S}\big[f\big](t,\zeta)\Big|^2 \dfrac{d\zeta\,dt}{\zeta^2}\right\}^{1/2} \left\{\int_{\mathbb R} \omega^2\left|\widehat f_{\bf\Lambda_{S}}(\omega)\right|^2 d\omega\right\}^{1/2}\ge \dfrac{ \sqrt C_{\psi}\,|B|}{2}\,\Big\|f\Big\|_2^{2}.\tag{4.7}
\end{align*}

\parindent=0mm \vspace{.1in}
{\it Proof.} The classical Heisenberg-Pauli-Weyl inequality in the SAFT domain is given by
\begin{align*}
\left\{\int_{\mathbb R} t^2\big| f(t)\big|^2 dt\right\}^{1/2}\left\{\int_{\mathbb R}\omega^2\left|{\mathcal F}_{SAFT}\big[f\big](\omega)\right|^2 d{\bf \omega}\right\}^{1/2}\ge \dfrac{|B|}{2} \left\{\int_{\mathbb R} \big| f(t)\big|^2 dt\right\}.\tag{4.8}
\end{align*}

Identifying ${\mathcal{A}}^{\psi}_{\bf\Lambda_S}\big[f\big](t,\zeta)$ as a function of the time variable $t$ and invoking (4.8), so that
\begin{align*}
&\left\{\int_{\mathbb R} t^2\big| {\mathcal{A}}^{\psi}_{{\bf\Lambda}_S}\big[f\big](t,\zeta)\big|^2 dt\right\}^{1/2}\\\
&\qquad\qquad\qquad\qquad\times\left\{\int_{\mathbb R}\omega^2\left|{\mathcal F}_{SAFT}\big[{\mathcal{A}}^{\psi}_{\bf\Lambda_S}\big[f\big](t,\zeta)\big](\omega)\right|^2 d{\omega}\right\}^{1/2}\\\\
&\qquad\qquad\qquad\qquad\qquad\qquad\qquad\qquad\qquad\ge \dfrac{|B|}{2} \left\{\int_{\mathbb R} \big| {\mathcal{A}}^{\psi}_{\bf\Lambda_S}\big[f\big](t,\zeta)\big|^2 dt\right\}.\tag{4.9}
\end{align*}
Integrating (4.9) with respect to the $d\zeta/\zeta^2$, we obtain
\begin{align*}
&\int_{\mathbb R^+} \left\{\int_{\mathbb R} t^2\big| {\mathcal{A}}^{\psi}_{\bf\Lambda_S}\big[f\big](t,\zeta)\big|^2 dt\right\}^{1/2}\\\
&\qquad\qquad\qquad\times\left\{\int_{\mathbb R}\omega^2\left|{\mathcal F}_{SAFT}\big[{\mathcal{A}}^{\psi}_{\Lambda_S}\big[f\big](t,\zeta)\big](\omega)\right|^2 d{\omega}\right\}^{1/2}\dfrac{d\zeta}{\zeta^2}\\\\
&\qquad\qquad\qquad\qquad\qquad\qquad\qquad\qquad\ge \dfrac{|B|}{2} \left\{\int_{\mathbb R\times\mathbb R^+} \big| {\mathcal{A}}^{\psi}_{\bf\Lambda_S}\big[f\big](t,\zeta)\big|^2 \dfrac{dt\,d\zeta}{\zeta^2}\right\}.\tag{4.10}
\end{align*}

By virtue of the Cauchy-Schwartz's inequality and Fubini theorem  we can express (4.10) as
\begin{align*}
&\left\{\int_{\mathbb R\times\mathbb R^+} t^2\big| {\mathcal{A}}^{\psi}_{\bf\Lambda_S}\big[f\big](t,\zeta)\big|^2 \dfrac{d\zeta\,dt}{\zeta^2}\right\}^{1/2}\\\
&\qquad\qquad\qquad\times\left\{\int_{\mathbb R\times\mathbb R^+} \omega^2\left|{\mathcal F}_{SAFT}\big[{\mathcal{A}}^{\psi}_{\bf\Lambda_S}\big[f\big](t,\zeta)\big](\omega)\right|^2 \dfrac{d\zeta\,d{\omega}}{\zeta^2}\right\}^{1/2}\\\\
&\qquad\qquad\qquad\qquad\qquad\qquad\qquad\qquad\qquad\qquad\qquad\qquad\quad\ge\dfrac{ C_{\psi}|B|}{2} \Big\|f\Big\|_2^2.
\end{align*}

\parindent=0mm \vspace{.1in}
Using Theorem 2.5, the above inequality  can be written as
\begin{align*}
&\Bigg\{\int_{\mathbb R\times\mathbb R^+} t^2\big| {\mathcal{A}}^{\psi}_{\bf\Lambda_S}\big[f\big](t,\zeta)\big|^2 \dfrac{d\zeta\,dt}{\zeta^2}\Bigg\}^{1/2}\\\
&\qquad\qquad\qquad\qquad\times\left\{\int_{\mathbb R} \omega^2\left|\widehat f_{\bf\Lambda_{S}}(\omega)\right|^2\left\{\int_{\mathbb R^+} \dfrac{\Big|\widehat \Psi_{\bf\Lambda_{S}}(\zeta\omega,\zeta)\Big|^2}{\zeta}\,d\zeta\right\}d{\omega}\right\}^{1/2}\\\\
&\qquad\qquad\qquad\qquad\qquad\qquad\qquad\quad\qquad\qquad\quad\qquad\qquad\ge\dfrac{C_{\psi}|B|}{2} \Big\|f\Big\|_2^2.
\end{align*}

Applying Theorem 2.5, we obtain the desired result:
\begin{align*}
\left\{\int_{\mathbb R\times\mathbb R^+} t^2\big| {\mathcal{A}}^{\psi}_{\bf\Lambda_S}\big[f\big](t,\zeta)\big|^2 \dfrac{d\zeta\,dt}{\zeta^2}\right\}^{1/2}\left\{\int_{\mathbb R} \omega^2\left|\widehat f_{\bf\Lambda_{S}}(\omega)\right|^2d{\omega}\right\}^{1/2}\ge\dfrac{\sqrt{C_{\psi}}\,|B|}{2} \Big\|f\Big\|_2^2.
\end{align*}

\parindent=0mm \vspace{.1in}
This completes the proof. \qquad\fbox

\parindent=0mm \vspace{.1in}

{\it Remark:} It should be noted that by choosing an appropriate matrix $\bf\Lambda_S$ yields the respective uncertainty inequalities for the various novel integral transforms.

\parindent=0mm \vspace{.1in}

\parindent=0mm \vspace{.2in}
{\bf{References}}

\end{document}